\documentclass{article}

\PassOptionsToPackage{numbers, sort&compress}{natbib}

\usepackage{natbib}
\usepackage[preprint]{neurips_2022}

\usepackage[utf8]{inputenc} 
\usepackage[T1]{fontenc}    
\usepackage{hyperref}       
\usepackage{url}            
\usepackage{booktabs}       
\usepackage{amsfonts}       
\usepackage{nicefrac}       
\usepackage{microtype}      
\usepackage{xcolor}         

\usepackage{times}
\usepackage{amsmath, nccmath, amssymb}
\usepackage{amsthm}
\usepackage{geometry}
\usepackage{comment}
\usepackage{enumitem}
\usepackage{color}
\usepackage{graphicx}
\usepackage{bbold}
\usepackage{xfrac}
\usepackage{thmtools}
\usepackage[ruled]{algorithm2e}
\usepackage{algorithmic}
\usepackage{mathtools}

\usepackage{notation}

\newtheorem{example}{Example} 
\newtheorem{theorem}{Theorem}
\newtheorem{lemma}{Lemma}
\newtheorem{remark}{Remark}

 \title{A gradient estimator via L1-randomization for online zero-order optimization with two point feedback}
 
\author{Arya Akhavan\\
      Istituto Italiano di Tecnologia\\
      CREST, ENSAE, IP Paris\\
      \texttt{aria.akhavanfoomani@iit.it}
      \And
      Evgenii Chzhen \\
      Université Paris-Saclay, CNRS\\
      Laboratoire de mathématiques d’Orsay\\
      \texttt{evgenii.chzhen@cnrs.fr}
      \AND
      Massimiliano Pontil \\
    Istituto Italiano di Tecnologia\\
      University College London\\
      \texttt{massimiliano.pontil@iit.it}
      \And
      Alexandre B. Tsybakov \\
      CREST, ENSAE, IP Paris \\
      \texttt{alexandre.tsybakov@ensae.fr}
      }

\begin{document}

\maketitle

\begin{abstract}
This work studies online zero-order optimization of convex and Lipschitz functions. We present  a novel gradient estimator based on two function evaluations and randomization on the $\ell_1$-sphere.
Considering different geometries of feasible sets and Lipschitz assumptions we analyse online dual averaging algorithm with our estimator in place of the usual gradient.
We consider two types of  assumptions on the noise of the zero-order oracle: canceling noise and adversarial noise.
We provide an anytime and completely data-driven algorithm, which is adaptive to all parameters of the problem.
In the case of canceling noise that was previously studied in the literature, our guarantees are either comparable or better than state-of-the-art bounds obtained by~\citet{duchi2015} and \citet{Shamir17} for non-adaptive algorithms. 
Our analysis is based on deriving a new 
weighted Poincaré type inequality for the uniform measure on the $\ell_1$-sphere with explicit constants, which may be of independent interest.
\end{abstract}

\section{Introduction}

In this work we study the problem of convex online zero-order optimization with two-point feedback, in which 
adversary fixes a sequence $f_1, f_2, \ldots : \bbR^d \to \bbR$ of convex functions and the goal of the learner is to minimize the cumulative regret with respect to the best action in a prescribed convex set $\com \subseteq \mathbb{R}^d$.  
This problem has received significant attention in the context of continuous bandits and online optimization \citep[see e.g.,][and references therein]{agarwal2010,flaxman2004,saha2011, bubeck2012regret, gasnikov_gradient-free_2016, Shamir17, bubeck2017kernel, akhavan2020, Gasnikov, lattimore21a}.

We consider the following protocol:
at each round $t = 1, 2, \ldots$ the algorithm chooses $\bx_t', \bx_t'' \in \bbR^d$ (that can be queried outside of $\com$) and the adversary reveals
    \begin{align*}        f_t(\bx_t') + \xi_t' \qquad\text{and}\qquad f_t(\bx_t'') + \xi_t''\enspace,
    \end{align*}
where $\xi'_t, \xi''_t \in \bbR$ are the noise variables (random or not) to be specified.    
Based on the above information and the previous rounds, the learner outputs $\bx_t \in \com$ and suffers loss $f_t(\bx_t)$.
The goal of the learner is to minimize the cumulative regret
\begin{align*}
    \sum_{t = 1}^T f_t(\bx_t) - \min_{\bx \in \com}\sum_{t = 1}^T f_t(\bx)\enspace.
\end{align*}
At the core of our approach is a novel zero-order gradient estimator
based on two function evaluations outlined in Algorithm \ref{algo:simplex}. A key novelty of our estimator is that 
it employs a randomization step over the $\ell_1$ sphere.
This is in contrast to most of the prior work \cite[see e.g.,][]{NY1983,PT90,flaxman2004,agarwal2010,Shamir13,duchi2015,BP2016,Gasnikov2017,akhavan2020,akhavan2021distributed} that was employing $\ell_2$ or $\ell_{\infty}$ type randomizations to define $\bx'_t, \bx''_t$.
We use the proposed estimator within an online dual averaging procedure to tackle the zero-order online convex optimization problem, matching or improving the state-of-the-art results.
\citet{duchi2015} and \citet{Shamir17} have studied instances of the above problem under the assumption that $\xi'_t = \xi''_t$, which we will further refer to as canceling noise assumption. Specifically, \cite{duchi2015} considered the stochastic optimization framework where $f_t= f$, for every $t$, and obtained bounds on the optimization error rather than on cumulative regret, while \cite{Shamir17} analyzed the case $\xi'_t = \xi''_t=0$.  The results in \cite{duchi2015, Shamir17} are obtained for the objective functions that are Lipschitz with respect to the $\ell_q$-norm for $q=1$ and $q=2$, although, with extra derivations it is possible to extend the above mentioned results beyond such cases.
The proposed method allows us to improve upon these results in several aspects.

{\bf Contributions.~} 
The contributions of the present paper can be summarized as follows. 
{\bf 1)} We present a new randomized zero-order gradient estimator and study its statistical properties, both under canceling noise and under adversarial noise (see Lemma~\ref{unb-est} and Lemma~\ref{var});
{\bf 2)} In the canceling noise case ($\xi'_t = \xi''_t$) in Theorem~\ref{no-noise-reg} we show that dual averaging based on our gradient estimator either improves or matches the state-of-the-art bounds~\cite{duchi2015,Shamir17}. We derive the results for Lipschitz functions with respect to all $\ell_q$-norms, $q \in [1,\infty]$. In particular, when $q=1$ and $\com$ is the probability simplex, our bound is better by a $\sqrt{\log(d)}$ factor than that of~\cite{duchi2015,Shamir17};
{\bf 3)} We propose a completely data-driven and anytime version of the algorithm, which is adaptive to all parameters of the problem. We show that it achieves analogous performance as the non-adaptive algorithm in the case of canceling noise and only slightly worse performance under adversarial noise. To the best of our knowledge, no adaptive algorithms were developed for zero-order online problems in our setting so far;
{\bf 4)} As a key element of our analysis, we derive in Lemma~\ref{lem:poincare_like_2} a weighted Poincaré type inequality~\citep[following the terminology of][]{BobkovLedoux09} with explicit constants for the uniform measure on $\ell_1$-sphere. This result may be of independent interest.

{\bf Notation.~} Throughout the paper we use the following notation.
We denote by $\norm{\cdot}_p$ the $\ell_p$-norm in $\bbR^d$. For any $\bx \in \bbR^d$ we denote by $\bx \mapsto \sign(\bx)$ the component-wise sign function (defined at $0$ as $1$). We let $\langle \cdot, \cdot \rangle$ be the standard inner product in $\mathbb{R}^d$.
For $p \in [1, \infty]$ we introduce the open $\ell_p$-ball 
and $\ell_p$-sphere respectively as
\begin{align*}
    \ball^d_p \triangleq \enscond{\bx \in \bbR^d}{\norm{\bx}_p < 1}\qquad\text{and}\qquad\sphere^d_p \triangleq \enscond{\bx \in \bbR^d}{\norm{\bx}_p = 1}\enspace.
\end{align*}
For two $a, b \in \bbR$, we denote by $a \wedge b$ (\emph{resp.} $a \vee b$) the minimum (\emph{resp.} the maximum) between $a$ and $b$. We denote by $\Gamma : (0, \infty) \to \bbR$, the gamma function. In what follows, $\log$ always stands for the natural logarithm and $e$ is Euler's number.

\section{The algorithm}

 Let $\Theta$ be a  closed convex subset of $\bbR^d$ and let $V : \Theta \to \bbR$ be a convex function.  
The  procedure that we propose in this paper is summarized in Algorithm \ref{algo:simplex}.

\begin{algorithm}[t!]
    \DontPrintSemicolon
   \caption{Zero-Order $\ell_1$-Randomized Online Dual Averaging}
  \label{algo:simplex}
   \SetKwInput{Input}{Input}
   \SetKwInOut{Output}{Output}
   \SetKwInput{Initialization}{Initialization}
   \Input{Convex function $V(\cdot)$, step size $\eta_1 > 0$, and parameters $h_t$, for $t=1, 2, \ldots,$}
   \Initialization{Generate independently vectors $\bzeta_1, \bzeta_2, \ldots$ uniformly distributed on $\sphere^d_1$, and set $\bz_1 = \bzero$}
   \For{ $t = 1, \ldots, $}{
   $\bx_t = \argmax_{\bx \in \com} \ens{\eta_t\scalar{\bz_t}{\bx} - V(\bx)}$\\
     \vspace{0.5mm}
   $y_t' = f_t(\bx_t+h_t\bzeta_t)+ \xi_t'\quad$ and $\quad y_t'' = f_t(\bx_t-h_t\bzeta_t)+ \xi_t''$\tcp*{Query}
   $\bg_t = \frac{d}{2h_t} (y_t'- y_t'') \sign(\bzeta_t) $\tcp*{$\ell_1$-gradient estimate}
   $\bz_{t+1}= \bz_t - \bg_t$\tcp*{Update $\bz_t$}
   update the step-size $\eta_{t+1}$}
\end{algorithm}

{\bf Intuition behind the gradient estimate.~}
The form of gradient estimator $\bg_t$ in Algorithm 1 is explained by Stokes' theorem (see Theorem~\ref{thm:ipp} in the appendix and the discussion that follows).
Stokes' theorem provides a connection between the gradient of a function $f$  (first order information) and $f$ itself (zero order information). Under some regularity conditions, it establishes that
\begin{align*}
        \int_{D} \nabla f(\bx) \d \bx = \int_{\partial D} f(\bx) \bn(\bx) \d S(\bx)\enspace,
    \end{align*}
    where $\partial D$ is the boundary of $D$, $\bn$ is the outward normal vector to $\partial D$, and $\d S(\bx)$ denotes the surface measure. Introducing $\bU^D$ and $\bzeta^{\partial D}$ distributed uniformly on $D$ and $\partial D$ respectively, we can rewrite the above identity as
    \begin{align*}
        \Exp[\nabla f(\bU^D)] = \frac{\vol_{d-1}(\partial D)}{\vol_d(D)} \cdot \Exp[ f(\bzeta^{\partial D}) \bn(\bzeta^{\partial D})]\enspace,
    \end{align*}
    where $\vol_{d-1}(\partial D)$ is the surface area of $D$ and $\vol_{d}(D)$ is its volume. In what follows we consider the special case  $D = \ball^d_1$. For this choice of $D$ we have $\bn(\bx) = \tfrac{1}{\sqrt{d}}\cdot\sign(\bx)$ with $\vol_{d-1}(\partial D) / \vol_{d}(D) = d^{3/2}$ leading to our gradient estimate for the two-point feedback setup.


{\bf Computational aspects.}
Let us highlight two appealing practical features of the $\ell_1$-randomized gradient estimator $\bg_t$ in Algorithm~\ref{algo:simplex}. First, 
we can easily evaluate any $\ell_p$-norm of $\bg_t$. Indeed, it holds that $\|\bg_t\|_{p} = (d^{1 + 1/p} /{2h_t}) |y_t' - y_t''|$, i.e., computing $\|\bg_t\|_{p}$ only requires $O(1)$ elementary operations. 
Second, this  gradient estimator is very economic in terms of the required memory: in order to store $\bg_t$ we only need $d$ bits and $1$ float. None of these properties is inherent to the popular alternatives based on the randomization over the $\ell_2$-sphere~\citep[see e.g.,][]{NY1983,flaxman2004,BP2016} or on  Gaussian randomization~\citep[see e.g.,][]{Nesterov2011,NS17,Ghadimi2013}.\\
To compute $\bg_t$ one needs to generate $\bzeta_t$ distributed uniformly on $\partial B_1^d$. The most straightforward way to do it consists in first generating a $d$-dimensional vector of i.i.d. centered scaled Laplace random variables and then normalizing this vector by its $\ell_1$-norm. The result is guaranteed to follow the uniform distribution on $\sphere_d^1$~\cite[see e.g.,][Lemma 1]{Schechtman_Zinn90}. Furthermore, to sample from the centered scaled Laplace distribution one can simply use inverse transform sampling. Indeed, if $U$ is distributed uniformly on $(0, 1)$, then $\log(2U)\ind{U > 1/2} - \log(2 - 2U)\ind{U \geq 1/2}$ follows centered scaled Laplace distribution.

\section{Assumptions}
\label{sec:ass_algo}
We say that the convex function $V(\cdot)$ is $1$-strongly convex with respect to the $\ell_p$-norm on $\Theta$ if 
\begin{align*}
    V(\bx') \geq V(\bx) + \scalar{\bw}{\bx' - \bx} + \frac{1}{2}\|\bx - \bx'\|_p^2\enspace,
\end{align*}
for all $\bx, \bx' \in \Theta$ and all $\bw \in \partial V(\bx)$, where $\partial V(\bx)$ is the subdifferential of $V$ at point $\bx$.

Throughout the paper, we assume that $p, q \in [1, \infty]$, $d \geq 3$, and set $p^*, q^* \in [1, \infty]$ such that $\tfrac{1}{p} + \tfrac{1}{p^*} = 1$ and $\tfrac{1}{q} + \tfrac{1}{q^*} = 1$, with the usual convention $1/\infty = 0$. We will use  the following assumptions.
\begin{assumption}\label{main-ass}
The following conditions hold:
\vspace{-2mm}
\begin{enumerate}
\itemsep0em
    \item The set $\Theta \subset \bbR^d$ is compact and convex.
    \item There exists $V : \com \to \bbR$, which is lower semi-continuous, $1$-strongly convex on $\com$ w.r.t. the $\ell_p$-norm and such that \[\sup_{\bx \in \com} V(\bx) - \inf_{\bx \in \com} V(\bx) \leq R^2\]
    for some constant $R>0$.
    \item Each function $f_t : \bbR^d \to \bbR$ is convex on $\bbR^d$ for all $t \geq 1$.
    \item For all $\bx, \bx' \in \bbR^d$, and all $t \geq 1$ we have $|f_t(\bx) - f_t(\bx')| \leq L \|\bx - \bx' \|_{q}$ for some constant $L>0$.
\end{enumerate}
\end{assumption}
\vspace{-2mm}
Assumption \ref{main-ass} is rather standard in the study of dual averaging-type algorithms and have been previously considered in the context of zero-order problems in~\cite{duchi2015,Shamir17}.
We assume that $\Theta$ is compact as we are interested in the worst-case regret, which ensures that $R < +\infty$. We discuss extensions of our results to the case of unbounded $\Theta$ in Section~\ref{sec:discission}.
Note that the constant $R > 0$ is {not} necessarily dimension independent.
Below we provide two classical examples of $V$~\cite[see e.g.,][Section 2]{Shalev-Shwartz}.

\begin{example}\label{example1}
Let $\com$ be any convex subset of $\bbR^d$ and $p \in (1, 2]$. Then, $V(\bx) = \tfrac{1}{2(p-1)}\|\bx\|_p^2$ is $1$-strongly convex on $\com$ w.r.t. the $\ell_p$-norm.
\end{example}
\begin{example}
\label{ex:simplex_setup}
Let $\com = \enscond{\bx \in \bbR^d}{\|\bx\|_1 = 1\,,\,\,\,\, \bx \geq 0}$. Then\footnote{We use the convention that $0\log(0) = 0$.}, $ V(\bx) = \sum_{j = 1}^d x_j \log(x_j)$ is $1$-strongly convex on $\com$ w.r.t. the $\ell_1$-norm and 
$R^2 \leq \log (d)$.
\end{example}

{\bf Assumptions on the noise.~} We consider two  different  assumptions on the noises $\xi'_t, \xi_t''$.
The first noise assumption is common in the stochastic optimization context~\cite[see e.g.,][]{duchi2015,Shamir17,Ghadimi2013,Nesterov2011,NS17}.
\begin{assumption}[Canceling noise]\label{ass:no_noise}
    For all $t = 1, 2, \ldots$, it holds that $\xi_t' = \xi_t''$ almost surely. 
\end{assumption}
\vspace{-2mm}
Formally, Assumption \ref{ass:no_noise} permits noisy evaluations of function values. However, due to the fact that we are allowed to query $f_t$ at two points, taking difference of $y_t'$ and $y_t''$ in the estimator of the gradient effectively erases the noise. It  results in a smaller variance of our
gradient estimator. Importantly, Assumption \ref{ass:no_noise} covers the case of no noise, that is, the classical online optimization setting as defined, e.g., in \cite{Shalev-Shwartz}.

Second, we consider an adversarial noise assumption, which is essentially equivalent to the assumptions used  in~\cite{akhavan2020,akhavan2021distributed}.  
\textbf{
\begin{assumption}[Adversarial noise]\label{ass:noise}
    For all $t = 1, 2, \ldots$, it holds that: (i) $\Exp[(\xi_t')^2] \leq \sigma^2$ and  $\Exp[(\xi_t'')^2] \leq \sigma^2$; (ii) $(\xi_t')_{t\ge1}$ and $(\xi_t'')_{t\ge 1}$ are independent of $(\bzeta_t)_{t\ge 1}$
\end{assumption}}
Assumption \ref{ass:noise} allows for stochastic $\xi_t'$ and $\xi_t''$ that are not necessarily zero-mean or independent over the trajectory. Furthermore, it permits bounded non-stochastic adversarial noises. 
Part (ii) of Assumption \ref{ass:noise} is always satisfied. Indeed, $\xi_t'$'s and $\xi_t''$'s are coming from the environment and are unknown to the learner while $\bzeta_t$'s are artificially generated by the learner. We mention part (ii) only for formal mathematical rigor.

Note that, since the choice of function $V$ belongs to the learner and  $\Theta$ is given, it is always reasonable to assume that parameter $R$ is known. At the same time, parameters $L$ and $\sigma$ may be either known or unknown.
We will study both cases in the next sections.

\section{Upper bounds on the regret}
\label{sec:rates}
In this section, we present the main convergence results for Algorithm~\ref{algo:simplex} when $L, \sigma, T$ are known to the learner. The case when they are unknown is analyzed in Section~\ref{adsect}, where we develop fully adaptive versions of Algorithm~\ref{algo:simplex}.

To state our results in a unified way, we introduce the following sequence that depends on the dimension $d$ and on the norm index $q \geq 1$: 
\begin{align*}
\mathrm{b}_{q}(d) \triangleq \frac{1}{d+1}\cdot\begin{dcases}
     q d^{\frac{1}{q}} &\text{if } q \in [1, \log(d)),\\
     e\log(d) &\text{if } q \geq \log(d).
     \end{dcases}
 \end{align*}    
The value $\mathrm{b}_{q}(d)$ will explicitly influence the choice of the step size $\eta > 0$ and of the discretization parameter $h > 0$.

The first result of this section establishes the convergence guarantees under the canceling noise assumption. This case was previously considered by~\citet{duchi2015} and \citet{Shamir17}. 
\begin{theorem}\label{no-noise-reg} Let Assumptions \ref{main-ass} and \ref{ass:no_noise} be satisfied. Then, Algorithm \ref{algo:simplex} with the parameters
\begin{align*}
    \eta = \frac{\cst   R}{L}\sqrt{\frac{d^{-1-\frac{2}{{q \wedge 2}}+\frac{2}{p}}}{T}}
    \quad\text{and any}\quad
    h\leq \frac{7R}{100\mathrm{b}_{q}(d)\sqrt{T}}d^{\frac{1}{2}+\frac{1}{q \wedge 2}-\frac{1}{p}}\enspace,
\end{align*}
where $\cst = (\sqrt{6} + \sqrt{12})^{-1}$,
satisfies, for any $\bx \in \com$,
\begin{align*}
    \Exp\left[\sum_{t=1}^{T}\big(f_{t}(\bx_t) - f_{t}(\bx)\big)\right] \leq 11.9\cdot RL \sqrt{Td^{1+\frac{2}{q \wedge 2}-\frac{2}{p}}} 
    \enspace.
\end{align*}
\end{theorem}
Note that, as in other related works~\cite{duchi2015,Shamir17,liang2014zeroth,Nesterov2011,NS17}, under the canceling noise (or no noise) assumption the discretization parameter $h > 0$ can be chosen arbitrary small. This is due to the fact that, under the canceling noise assumption, the variance of the gradient estimate $\bg_t$ is bounded by a constant independent of $h$. It is no longer the case under the adversarial noise assumption as exhibited in the next theorem.

\begin{theorem}\label{noisy-reg} Let Assumptions \ref{main-ass} and \ref{ass:noise} be satisfied. Then Algorithm \ref{algo:simplex} with the parameters
\begin{align*}
    \eta = \frac{R}{\sqrt{TL}}\parent{\frac{\sigma\mathrm{b}_{q}(d)}{\sqrt{2}R} \sqrt{Td^{4-\frac{2}{p}}}{+}\cst  Ld^{1+\frac{2}{{q \wedge 2}}-\frac{2}{p}}}^{-\frac{1}{2}}
    \quad\text{and}\quad
    h= \left(\frac{\sqrt{2}R\sigma}{L\mathrm{b}_{q}(d)}\right)^{\frac{1}{2}}T^{-\frac{1}{4}}d^{1-\frac{1}{2p}}\enspace,
\end{align*}
where $\cst = 6(1{+}\sqrt{2})^2$, 
satisfies, for any $\bx \in \com$,
\begin{align*}
   \Exp\left[\sum_{t=1}^{T}\big(f_{t}(\bx_t) - f_{t}(\bx)\big)\right] \leq 11.9&\cdot RL\sqrt{Td^{1+\frac{2}{{q \wedge 2}}-\frac{2}{p}}}\\
   +
   &\,\,2.4\cdot\sqrt{RL\sigma}T^{\frac{3}{4}}\cdot \begin{dcases}
     \sqrt{q d^{1+\frac{1}{q}-\frac{1}{p}}} &\text{if } q \in [1, \log(d)),\\
     \sqrt{e\log(d)d^{1-\frac{1}{p}}} &\text{if } q \geq \log(d).
     \end{dcases}
\end{align*}
\end{theorem}

{\bf Comparison to~~state-of-the-art bounds.~}
We provide two examples of $p, q, \com$ and compare results for our new method to those of~\cite{duchi2015,Shamir17} where only  the canceling noise Assumption~\ref{ass:no_noise} and $q\in\{1,2\}$ were considered.
\begin{corollary}\label{cor1}
Let $p = q = 2$ and $\com = \ball^d_2$. Then under Assumption~\ref{ass:no_noise}, Algorithm~\ref{algo:simplex} with $V : \com \to \bbR$ defined in Example~\ref{example1}, satisfies
\begin{align*}
    \Exp\left[\sum_{t=1}^{T}\big(f_{t}(\bx_t) - f_{t}(\bx)\big)\right] \leq 11.9\cdot L \sqrt{dT} \enspace.
\end{align*}
\end{corollary}
In the setup of Corollary~\ref{cor1}, \citet{duchi2015} obtain $O(L \sqrt{dT\log(d)})$ rate and~\citet{Shamir17} exhibits~$O(L \sqrt{dT})$, which is the optimal rate. Both results do not specify the leading absolute constants.
\begin{corollary}\label{cor2}
Let $p = q = 1$ and $\com = \enscond{\bx \in \bbR^d}{\bx \geq 0,\,\, \|\bx\|_1 = 1}$. Then under Assumption~\ref{ass:no_noise}, Algorithm~\ref{algo:simplex} with $V : \com \to \bbR$ defined in Example~\ref{ex:simplex_setup}, satisfies
\begin{align*}
    \Exp\left[\sum_{t=1}^{T}\big(f_{t}(\bx_t) - f_{t}(\bx)\big)\right] \leq 11.9\cdot L \sqrt{dT\log(d)} \enspace.
\end{align*}
\end{corollary}
In the setup of Corollary~\ref{cor2}, \citet{Shamir17} proves the rate~$O(L \sqrt{dT}\log(d))$ for the method with $\ell_2$-randomization.  On the other hand,~\citet{duchi2015} derived a lower bound $\Omega(\sqrt{dT / \log(d)})$. Thus, our algorithm further reduces the gap between the upper and the lower bounds. 

Finally, note that in  the case $p=1$, $q=2$ with $V : \com \to \bbR$ defined in Example~\ref{ex:simplex_setup} the bound of Theorem \ref{no-noise-reg} is of the order~$O( \sqrt{T\log(d)})$. This case was 
handled by an algorithm with $\ell_1$-randomization slightly different from ours in 
\cite{gasnikov_gradient-free_2016} leading to the suboptimal rate~$O(\sqrt{dT\log(d)})$.

\section{Adaptive algorithms}\label{adsect}
Theorems~\ref{no-noise-reg} and~\ref{noisy-reg} used the step size $\eta$ and the discretization parameter $h$ that depend on the potentially unknown quantities $L$, $\sigma$, and the optimization horizon $T$. In this section, we show that, under the canceling noise Assumption~\ref{ass:no_noise}, adaptation to unknown $L$ comes with nearly no price. 
On the other hand, under the adversarial noise Assumption~\ref{ass:noise}, our adaptive rate has a slightly worse dependence on $L$ and $\sigma$ in the dominant term. The proof is based on combining the adaptive scheme for online dual averaging~\cite[see Section 7.13 in][for an overview]{Orabona2019} with our bias and variance evaluations, cf. Section \ref{sec:proofs} below.
\begin{theorem}\label{no-noise-ad} Let Assumptions \ref{main-ass} and \ref{ass:no_noise} be satisfied. Then, Algorithm \ref{algo:simplex} with the parameters\footnote{We adopt the convention that $\eta_1 = 1$ and $1 / 0 = 1$ in the definition of $\eta_t$.}
\begin{align*}
    \eta_t =  \frac{R}{\sqrt{2.75\cdot\sum_{k=1}^{t-1}\norm{\bg_k}^2_{p^*}}}
    \quad\text{and any}\quad
    h_t\leq\frac{7R}{200\mathrm{b}_{q}(d)\sqrt{t}}d^{\frac{1}{2}+\frac{1}{q \wedge 2}-\frac{1}{p}}\enspace,
\end{align*}
satisfies for any $\bx \in \com$
\begin{align*}
    \Exp\left[\sum_{t=1}^{T}\big(f_{t}(\bx_t) - f_{t}(\bx)\big)\right] \leq 110.6\cdot RL \sqrt{Td^{1+\frac{2}{q \wedge 2}-\frac{2}{p}}}
    \enspace.
\end{align*}

\end{theorem}
The above result gives, up to an absolute constant, the same convergence rate as that of the non-adaptive Theorem~\ref{no-noise-reg}. In other words, the price for adaptive algorithm does not depend on the parameters of the problem.
Finally, we derive an adaptive algorithm under Assumption~\ref{ass:noise}.
\begin{theorem}\label{noisy-reg-ad} Let Assumptions \ref{main-ass} and \ref{ass:noise} be satisfied. Then, Algorithm \ref{algo:simplex} with the parameters
\begin{align*}
    \eta_t =  \frac{R}{\sqrt{2.75\cdot\sum_{k=1}^{t-1}\norm{\bg_k}^2_{p^*}}}
    \quad\text{and any}\quad
    h_t = \left(6.65\sqrt{6}\cdot\frac{R}{\mathrm{b}_{q}(d)}\right)^{\frac{1}{2}}t^{-\frac{1}{4}}d^{1-\frac{1}{2p}}\enspace,
\end{align*}
satisfies for any $\bx \in \com$
\begin{align*}
   \Exp\left[\sum_{t=1}^{T}\big(f_{t}(\bx_t) - f_{t}(\bx)\big)\right] \leq 110.6&\cdot RL\sqrt{Td^{1+\frac{2}{{q \wedge 2}}-\frac{2}{p}}} \\
   &+
   5.9\cdot\sqrt{R}\left(\sigma {+} L\right)T^{\frac{3}{4}}\cdot \begin{dcases}
     \sqrt{q d^{1+\frac{1}{q}-\frac{1}{p}}} &\text{if } q \in [1, \log(d))\\
     \sqrt{e\log(d)d^{1-\frac{1}{p}}} &\text{if } q \geq \log(d)
     \end{dcases}\enspace.
\end{align*}
\end{theorem}
Note that the bound of Theorem \ref{noisy-reg-ad} has a less advantageous dependency on $\sigma$ and $L$ compared to Theorem~\ref{noisy-reg}, where we had $\sqrt{\sigma L}$ instead of $\sigma + L$. We remark that if $\sigma$ is known but $L$ is unknown, one can recover the $\sqrt{\sigma L}$ dependency by selecting $h_t$ depending on $\sigma$. We do not state this result that can be derived in a similar way and favor here only the fully adaptive version.

\section{Elements of proofs}
\label{sec:proofs}

In this section, we outline major ingredients for the proofs of  Theorems~\ref{no-noise-reg} -- \ref{noisy-reg-ad}. The full proofs can be found in Appendix~\ref{app:upper_bound}. Here, we only focus on novel elements without reproducing the general scheme of online dual averaging analysis~\cite[see e.g.,][]{Shalev-Shwartz, Orabona2019}. 
Namely, we highlight two key facts, which are the smoothing lemma (Lemma~\ref{unb-est}) and the weighted Poincaré type inequality for the uniform measure on $\sphere_1^d$  (Lemma~\ref{lem:poincare_like_2}) used to control the variance.


\subsection{Bias and smoothing lemma}
\label{sec:bias_smoothing_lemma}
First, as in the prior work that was using smoothing ideas \cite[see e.g.,][]{NY1983, flaxman2004, Shamir17}, we show that our gradient estimate $\bg_t$ is an unbiased estimator of a surrogate version of $f_t$ and establish its approximation properties.
\begin{lemma}[Smoothing lemma]\label{unb-est}
Fix $h > 0$ and $q \in [1,\infty]$.
Let $f: \bbR^d \to \bbR$ be an $L$-Lipschitz function w.r.t. the $\ell_{q}$-norm.
Let $\bU$ be distributed uniformly on $\ball^d_1$ and $\bzeta$ be distributed uniformly on $\sphere_d^1$.
Let ${\sf f}_{h}(\bx) \eqdef \Exp[f(\bx + h\bU)]$ for $\bx \in \bbR^d$. Then ${\sf f}_{h}$ is differentiable and
\begin{align*}
    \Exp\left[ \frac{d}{2h}\big(f(\bx + h \bzeta) - f(\bx - h \bzeta)\big) \sign(\bzeta)\right] = \nabla {\sf f}_{h}(\bx)\enspace.
\end{align*}
Furthermore, we have for all $d \geq 3$ and all $\bx \in \bbR^d$,
\begin{align}\label{eq:lem1}
     |{\sf f}_{h}(\bx) - f(\bx)|
     \leq
     \mathrm{b}_{q}(d) Lh\enspace.
\end{align}
Finally, if $\com \subset \bbR^d$ is convex, $f$ is convex in $\com {+} h  \ball_1^d$, then ${\sf f}_{h}$ is convex in $\com$ and ${\sf f}_h(\bx) \geq f(\bx)$ for $\bx \in \com$.

\end{lemma}
\begin{proof}
There are three claims to prove.
 For the first one, we notice that $\bzeta$ has the same distribution as $-\bzeta$, hence,
 \begin{align*}
     \Exp\left[ \frac{d}{2h}\big(f(\bx + h \bzeta) - f(\bx - h \bzeta)\big) \sign(\bzeta) \right] = \Exp\left[ \frac{d}{h}f(\bx + h \bzeta)\sign(\bzeta) \right]\enspace,
 \end{align*}
and the first claim follows from Theorem~\ref{thm:ipp2} in the Appendix (a version of Stokes', or divergence, theorem)  applied to $g(\cdot) = f(\bx + h \cdot)$ with observation that $\nabla g(\cdot) = h \nabla f(\bx + h \cdot)$ where $\nabla f$ is the gradient defined almost everywhere and whose existence is ensured by the Rademacher theorem.

We now prove the approximation property \eqref{eq:lem1}.  Assuming $d \geq 3$ grants that $\log(d) \geq 1$.
Since $f$ is $L$-Lipschitz w.r.t. the $\ell_{q}$-norm we get that,
for any $\bx \in \bbR^d$, 
\begin{align}
    \label{eq:unb-est1}
    |{\sf f}_{h}(\bx) - f(\bx)| \leq Lh\, \Exp \|\bU\|_q\enspace.
\end{align}
If $q \in [1, \log(d))$ then \eqref{eq:lem1} follows from Lemma~\ref{eq:moments_uniform_l1}. If $q \geq \log(d)$ then using again Lemma~\ref{eq:moments_uniform_l1} we find
\begin{align*}
    \Exp \|\bU\|_q \leq \Exp\|\bU\|_{\log(d)} \leq \frac{\log(d)d^{\frac{1}{\log(d)}}}{d + 1} = \frac{e \log(d)}{d + 1}\enspace,
\end{align*}
which together with~\eqref{eq:unb-est1} yields the desired bound.

Finally, if $f$ is convex in $\com + h \ball_1^d$, then for all $\bx, \bx' \in \com$ and $\alpha \in [0, 1]$ we have
\begin{align*}   {\sf f}_h(\alpha \bx + (1 - \alpha) \bx') =  \Exp\bigg[ f\big(\alpha (\bx + h\bU) + (1 - \alpha) (\bx' + h\bU ) \big)\bigg]  \leq \alpha {\sf f}_h(\bx) + (1 - \alpha) {\sf f}_h(\bx')\enspace.\end{align*}
Thus ${\sf f}_h$ is indeed convex on $\com$. Furthermore, again by convexity of $f$, we deduce that for any $\bx \in \com$
\begin{align*}
    {\sf f}_h(\bx) = \Exp[f(\bx + h\bU)] \geq \Exp\left[f(\bx) + \scalar{\bw}{h \bU} \right] = f(\bx)\quad\text{where $\bw \in \partial f(\bx)$}\,\,.
    &\qedhere
\end{align*}
\end{proof}
The proof of Lemma~\ref{unb-est} relies on the control of the $\ell_q$-norm of random vector $\bU$ established in the next result. 
\begin{lemma}
\label{eq:moments_uniform_l1}
Let $q \in [1,\infty)$ and let $\bU$ be distributed uniformly on $\ball^d_1$. Then $\Exp \|\bU\|_q \leq \tfrac{q d^{\frac{1}{q}}}{d+1}$.
\end{lemma}

\begin{proof}
Let $W_1, \ldots, W_d, W_{d+1}$ be \iid random variables having the Laplace distribution with mean $0$ and scale parameter $1$. Set $\bW = (W_1, \ldots, W_d)$.
Then, following \cite[Theorem 1]{Barthe_Guedon_Mendelson_Naor05} we have
\begin{align*}
    \bU\stackrel{d}{=} \frac{\bW}{\norm{\bW}_1 + |W_{d+1}|}\enspace,
\end{align*}
where $\stackrel{d}{=}$ denotes equality in distribution. Furthermore, \cite[Lemma 1]{Schechtman_Zinn90} states that the random variables
\begin{align*}   \frac{(\bW, |W_{d+1}|)}{\norm{\bW}_1 + |W_{d+1}|}\qquad\text{and}\qquad \norm{\bW}_1 + |W_{d+1}|\enspace, \end{align*}
are independent.
Hence, for any $q \in [1, \infty)$, it holds that
\begin{align*} \Exp \|\bU\|_q = \frac{\Exp \|\bW\|_q}{\Exp \|(\bW, W_{d+1})\|_1} = \frac{1}{d+1}\Exp \|\bW\|_q \stackrel{(a)} {\leq} \frac{1}{d+1}\bigg( \Exp \|\bW\|_q^q \bigg)^{\frac 1 q} = \frac{d^{\frac 1 q}\Gamma^{\frac 1 q}(q + 1)}{d+1} \stackrel{(b)}{\leq} \frac{qd^{\frac 1 q}}{d+1}\enspace, \end{align*}
where $(a)$ follows from Jensen's inequality and $(b)$ uses the fact that $\Gamma^{1/q}(q {+} 1) \leq q$ for $q \geq 1$.
\end{proof}

\subsection{Variance and weighted Poincaré type inequality}
\label{subsec:var_poincare}

We additionally need to control the squared $\ell_{p^*}$-norm of each gradient estimator $\bg_t$. This is where we get the main improvement of our procedure compared to previously proposed methods.
To derive the result, we first establish the following lemma of independent interest, which allows us to control the variance of Lipschitz functions on $\sphere^d_1$. The proof of this lemma is given in the Appendix.
\begin{lemma}
    \label{lem:poincare_like_2}
    Let $d \geq 3$. Assume that $G : \bbR^d \to \bbR$  is a continuously differentiable function, and $\bzeta$ is distributed uniformly on $\sphere^d_1$. Then
    \begin{align*}
        \Var(G(\bzeta)) \leq \frac{4}{d(d-2)}\Exp\left[\|\nabla G(\bzeta)\|_2^2\parent{1 + \sqrt{d}\|\bzeta\|_2}^2\right]\enspace.
    \end{align*}
     Furthermore, if $G : \bbR^d \to \bbR$ is an $L$-Lipschitz function w.r.t. the $\ell_2$-norm then 
    \begin{align*}
        \Var(G(\bzeta)) \leq \frac{4L^2}{d(d-2)}\parent{1 + \sqrt{\frac{2d}{d+1}}}^2\enspace.
    \end{align*}
\end{lemma}
\begin{remark}
Since $d^2 / (d(d-2)) \leq 3$ for all $d \geq 3$, the last inequality of Lemma \ref{lem:poincare_like_2} implies that
    \begin{align}
    \label{eq:bound_to_use_for_variance}
        \Var(G(\bzeta)) \leq 12\big({1 + \sqrt{2}}\big)^2\big({L}/{d}\big)^2\,\,, \qquad \forall d \geq 3\enspace.
    \end{align}
\end{remark}
We can now deduce the following bound on the squared $\ell_{p^*}$-norm of $\bg_t$.
\begin{lemma}\label{var}
Let $p \in [1, \infty]$ and $p^* = \tfrac{p}{p-1}$.
Assume that $f_t$ is $L$-Lipschitz w.r.t. the $\ell_q$-norm. Then, for all $d \geq 3$, 
\begin{align*}
    \Exp\|\bg_t\|_{p^*}^2 \leq 12(1+\sqrt{2})^2L^2 d^{1 + \frac{2}{q \wedge 2} - \frac{2}{p}}
    +
    \begin{dcases}
    0 & \text{under canceling noise Assumption~\ref{ass:no_noise}},\\
    \frac{d^{4-\frac{2}{p}}\sigma^2}{h^2}& \text{under adversarial noise Assumption~\ref{ass:noise}}.
    \end{dcases}
\end{align*}
\end{lemma}
\begin{proof}
Using the definition of $\bg_t$ we get
\begin{align*}
    \Exp[\|\bg_t\|_{p^*}^2 \mid \bx_t]
    &= \frac{d^2}{4h^2}\Exp[(f_t(\bx_t + h \bzeta_t) - f_t(\bx_t - h \bzeta_t) + \xi_t' - \xi_t'')^2\|\sign(\bzeta_t)\|_{p^*}^2 \mid \bx_t]
    \\
    &=\frac{d^{4-\frac{2}{p}}}{4h^2}\Exp[(f_t(\bx_t + h \bzeta_t) - f_t(\bx_t - h\bzeta_t) + \xi_t' - \xi_t'')^2 \mid \bx_t]\enspace.
\end{align*}
Let $G(\bzeta) \triangleq f_t(\bx_t + h \bzeta) - f_t(\bx_t - h \bzeta)$. First, observe that $\Exp[G(\bzeta_t) \mid \bx_t] = 0$ and under both Assumption~\ref{ass:no_noise} and Assumption~\ref{ass:noise}(ii) it holds that $\Exp[G(\bzeta_t)(\xi_t' - \xi_t'') \mid \bx_t ] = 0$. Using these remarks and the fact that under adversarial noise Assumption~\ref{ass:noise}, $\Exp[(\xi_t' - \xi_t'')^2 \mid \bx_t] \leq 4\sigma^2$, we find:
\begin{align*}
    \Exp[\|\bg_t\|_{p^*}^2 \mid \bx_t] \leq \frac{d^{4-\frac{2}{p}}}{4h^2}  \left(\Var(G(\bzeta_t) \mid \bx_t) + \begin{dcases}
    0 & \text{under cancelling noise Assumption~\ref{ass:no_noise}}\\
    4\sigma^2 & \text{under adversarial noise Assumption~\ref{ass:noise}}
    \end{dcases}\right).
\end{align*}
Furthermore, since $f_t$ is $L$-Lipschitz, w.r.t. the $\ell_q$-norm, the map $\bzeta \mapsto G(\bzeta)$ is $\big(2Lhd^{\frac{1}{q \wedge 2} - \frac{1}{2}}\big)$-Lipschitz w.r.t. the $\ell_2$-norm. Applying~\eqref{eq:bound_to_use_for_variance} to bound $\Var(G(\bzeta_t) \mid \bx_t)$, yields the desired result.
\end{proof}
Note that under adversarial noise Assumption~\ref{ass:noise}, the bound on squared $\ell_{p^*}$-norm of $\bg_t$ gets an additional term $d^{4-\frac{2}{p}}\sigma^2h^{-2}$. In contrast to the case of canceling noise Assumption~\ref{ass:no_noise}, this does not allow us to take $h$ arbitrary small hence inducing the bias-variance trade-off.

\section{Numerical illustration}
\label{sec:exps}
\begin{figure}[t]
\includegraphics[width=0.5\textwidth]{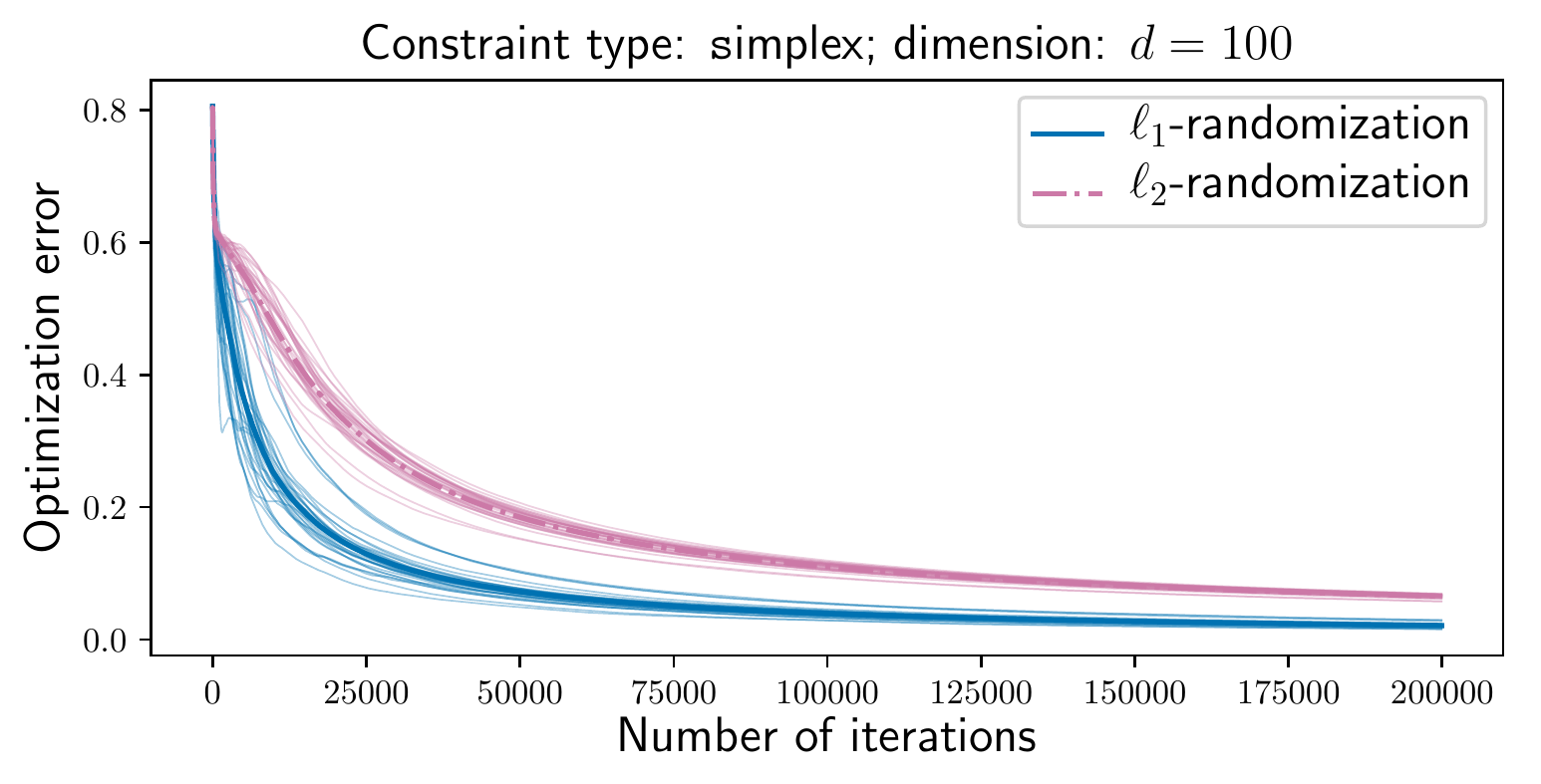}\hfill
\includegraphics[width=0.5\textwidth]{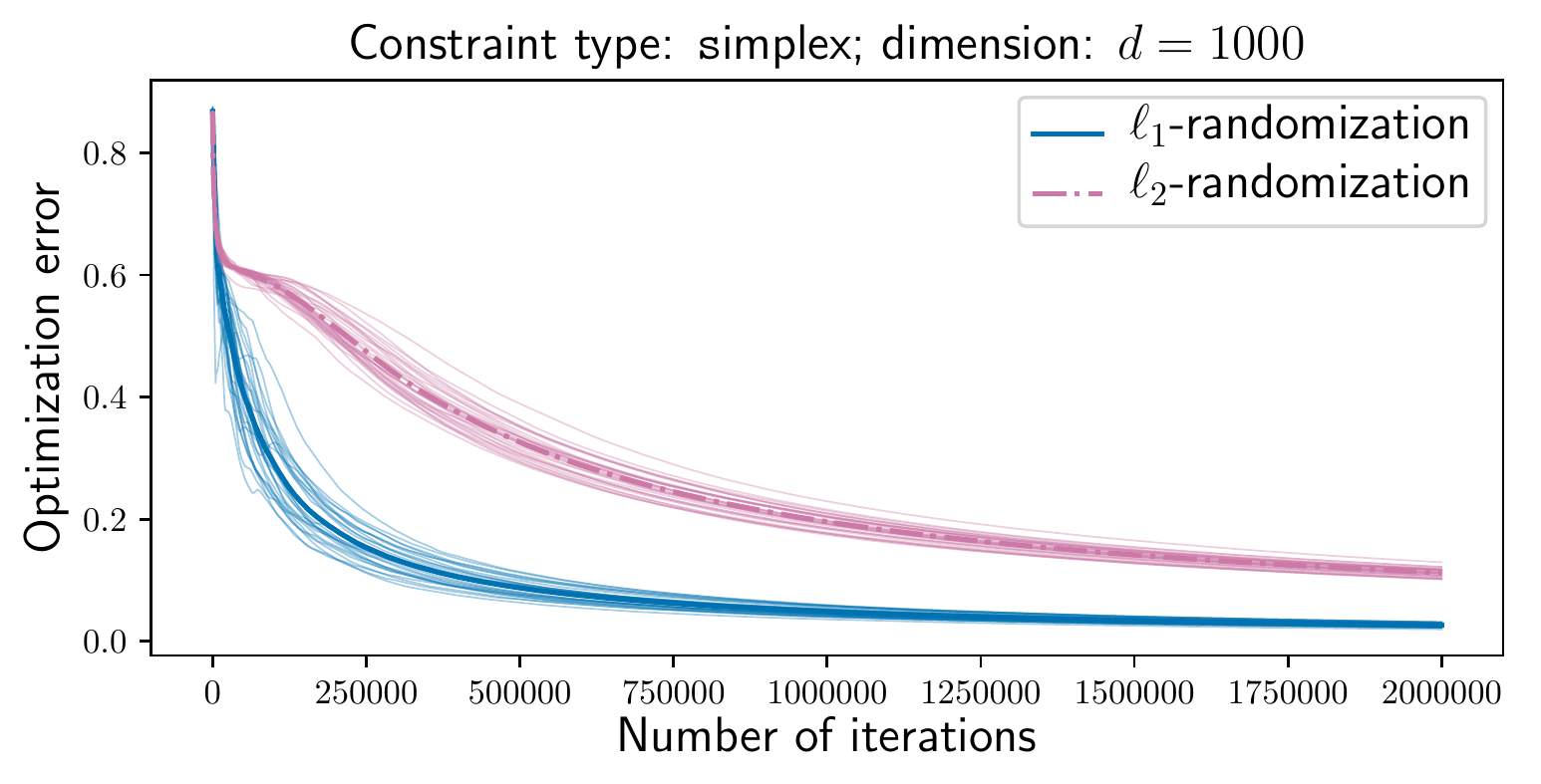}
\caption{Opt. error vs. number of iterations for $\ell_2$-randomization (as in~\cite{Shamir17}) and our method.
}
\label{fig:plots}
\end{figure}

In this section, we provide a numerical comparison of our algorithm with the method based on $\ell_2$-randomization from~\citet{Shamir17} (see Appendix~\ref{app:l2_description} for the definition). 
We consider the no noise model and $f_t  = f$, $\forall t$, with the function $f:\mathbb{R}^d \to \mathbb{R}$ defined as
\begin{align*}
    f(\bx) = \norm{\bx - \bc}_{2} + \norm{\bx - 0.1\cdot\bc}_{1}\enspace,
\end{align*}
where $\bc = (c_1, \ldots, c_d)^\top \in \bbR^d$ such that $c_j = {\exp(j)}/{\sum_{i = 1}^d \exp(i)}$ for $j = 1, \ldots, d$.
We choose
\[
\com = \enscond{\bx \in \bbR^d}{\|\bx\|_1 = 1,\, \bx \geq 0}~~~{\rm and}~~~V(\bx)=\sum_{j=1}^{d}x_j\log(x_j)\enspace.
\] 
As stated in Example \ref{ex:simplex_setup}, $V$ is $1$-strongly convex on $\com$ w.r.t. the $\ell_1$-norm and $R \leq \sqrt{\log(d)}$. Moreover, $f$ is a Lipschitz function w.r.t. the $\ell_1$-norm. 
We deploy the adaptive parameterization that appears in Theorem \ref{no-noise-ad}.
In Figure \ref{fig:plots} we present the optimization error of the algorithms, which is defined as \[
f\bigg(\frac{1}{t}\sum_{i = 1}^t\bx_i\bigg) - \min_{\bx \in \com}f(\bx)\enspace.
\]

The results are reported over $30$ trials. We plot all the $30$ runs alongside the average performance. One can observe that the $\ell_1$-randomization method behaves significantly better than the $\ell_2$-randomization algorithm.  The theoretical bound for our method in this setup has a $\sqrt{\log d}$ gain in the rate. 

\section{Discussion and comparison to prior work}
\label{sec:discission}
We introduced and analyzed a novel estimator for the gradient based on randomization over the $\ell_1$-sphere. We established guarantees for the online dual averaging algorithm with the gradient replaced by the proposed estimator. We provided an anytime and completely data-driven algorithm, which is adaptive to all parameters of the problem.
Our analysis is based on deriving a weighted
Poincaré type inequality for the uniform measure on the $\ell_1$-sphere that may be of independent interest.\\ 
Under the {\it canceling noise assumption} and $q\in\{1,2\}$, our setting is analogous to \cite{duchi2015, Shamir17}. For the case $q = p =2$ and canceling noise, we show that the performance of our method is the same as in \cite[Corollary 2]{Shamir17} up to absolute constants that were not made explicit in \cite{Shamir17}.
For the case of $q=p=1$ and canceling noise, we improved the bound \cite[Corollary 3]{Shamir17} by a $\sqrt{\log(d)}$ factor. For the case $q = 2$, $p \geq 1$, comparing with the lower bound in \cite[Proposition 1]{duchi2015}, shows that the result of Theorem \ref{no-noise-reg} is minimax optimal. For the case $q = p = 1$, \cite[Proposition 2]{duchi2015} shows that our result in Theorem \ref{no-noise-reg} is optimal up to a $\log(d)$ factor.\\
Under the {\it adversarial noise assumption}, Theorem \ref{noisy-reg} provides the rate $O(T^{3/4})$, that is, we get an additional $T^{1/4}$ factor compared to the canceling noise case. It remains unclear whether it is optimal under adversarial noise -- this question deserves further investigation. Note that, under sub-Gaussian \iid noise assumption and $q =p = 2$, one can achieve the rate $\tilde{O}(d^a\sqrt{T})$ with a relatively big $a>0$~\cite{agarwal2011,belloni2015escaping, bubeck2017kernel, lattimore21a}.
In particular, with an ellipsoid type 
method \cite{lattimore21a} obtains the rate $\mathcal{O}(d^{4.5}\sqrt{T}\log(T)^2)$ for the cumulative regret.\\
Finally, let us discuss the compactness of $\Theta$. 
It is straightforward to extend the results of Theorems~\ref{no-noise-reg},~\ref{noisy-reg} to any closed convex $\Theta$ considering the regret against a fixed action $\bx \in \Theta$. Indeed, using~\citep[Corollary 7.9]{Orabona2019}, one only needs to replace $R$ appearing in both Theorems~\ref{no-noise-reg},~\ref{noisy-reg} by an upper bound on $\sqrt{V(\bx) - \inf_{\bx' \in \Theta} V(\bx')}$. The adaptive case is more complicated. One way to tackle this case is to use~\citep[Theorem 1]{orabonapal16} requiring a control of $\Exp\max_{t =1, \ldots, T} \|\bg_t\|_{p^*}$. This term can be controlled under the canceling noise Assumption~\ref{ass:no_noise} using the Lipschitzness of $f_t$'s, so that Theorem~\ref{no-noise-ad} extends to unbounded $\Theta$. However, without the canceling noise assumption, following the approach outlined above, 
one needs to control $\Exp \max_{t = 1, \ldots, T} \tfrac{|\xi_t' - \xi''_{t}|}{h_t}$. The adversarial noise Assumption~\ref{ass:noise} is not sufficient to reasonably control this term, so that extending Theorem~\ref{noisy-reg-ad} to unbounded $\Theta$ is not possible without further assumptions. 

\begin{ack}
The work of E. Chzhen was supported by ANR PIA funding: ANR-20-IDEES-0002. The work of M.Pontil is partially supported by European Union's Horizon 2020 research and innovation programme under ELISE grant agreement No 951847. The research of A.Akhavan and A.B. Tsybakov is supported by a grant of the French National Research Agency (ANR), “Investissements d’Avenir” (LabEx Ecodec/ANR-11-LABX-0047).
\end{ack}

\bibliography{biblio}

\begin{thebibliography}{}

\bibitem[Agarwal et~al., 2010]{agarwal2010}
Agarwal, A., Dekel, O., and Xiao, L. (2010).
\newblock Optimal algorithms for online convex optimization with multi-point
  bandit feedback.
\newblock In {\em Proc. 23rd International Conference on Learning Theory},
  pages 28--40.

\bibitem[Agarwal et~al., 2011]{agarwal2011}
Agarwal, A., Foster, D.~P., Hsu, D.~J., Kakade, S.~M., and Rakhlin, A. (2011).
\newblock Stochastic convex optimization with bandit feedback.
\newblock In {\em Advances in Neural Information Processing Systems},
  volume~25, pages 1035--1043.

\bibitem[Akhavan et~al., 2020]{akhavan2020}
Akhavan, A., Pontil, M., and Tsybakov, A. (2020).
\newblock Exploiting higher order smoothness in derivative-free optimization
  and continuous bandits.
\newblock In {\em Advances in Neural Information Processing Systems},
  volume~33.

\bibitem[Akhavan et~al., 2021]{akhavan2021distributed}
Akhavan, A., Pontil, M., and Tsybakov, A.~B. (2021).
\newblock Distributed zero-order optimization under adversarial noise.
\newblock In {\em Advances in Neural Information Processing Systems},
  volume~34.

\bibitem[Bach and Perchet, 2016]{BP2016}
Bach, F. and Perchet, V. (2016).
\newblock Highly-smooth zero-th order online optimization.
\newblock In {\em Proc. 29th Annual Conference on Learning Theory}, pages
  1--27.

\bibitem[Barthe et~al., 2005]{Barthe_Guedon_Mendelson_Naor05}
Barthe, F., Gu{\'e}don, O., Mendelson, S., and Naor, A. (2005).
\newblock {A probabilistic approach to the geometry of the Lpn-ball}.
\newblock {\em The Annals of Probability}, 33(2):480 -- 513.

\bibitem[Barthe and Wolff, 2009]{barthe2009remarks}
Barthe, F. and Wolff, P. (2009).
\newblock Remarks on non-interacting conservative spin systems: the case of
  gamma distributions.
\newblock {\em Stochastic processes and their applications}, 119(8):2711--2723.

\bibitem[Belloni et~al., 2015]{belloni2015escaping}
Belloni, A., Liang, T., Narayanan, H., and Rakhlin, A. (2015).
\newblock Escaping the local minima via simulated annealing: Optimization of
  approximately convex functions.
\newblock In {\em Proc. 28th Annual Conference on Learning Theory}, pages
  240--265.

\bibitem[Bobkov and Ledoux, 1997]{bobkov1997poincare}
Bobkov, S. and Ledoux, M. (1997).
\newblock Poincar{\'e}’s inequalities and {T}alagrand’s concentration
  phenomenon for the exponential distribution.
\newblock {\em Probability Theory and Related Fields}, 107(3):383--400.

\bibitem[Bobkov and Ledoux, 2009]{BobkovLedoux09}
Bobkov, S.~G. and Ledoux, M. (2009).
\newblock {Weighted Poincaré-type inequalities for Cauchy and other convex
  measures}.
\newblock {\em The Annals of Probability}, 37(2):403 -- 427.

\bibitem[Boucheron et~al., 2013]{boucheron2013concentration}
Boucheron, S., Lugosi, G., and Massart, P. (2013).
\newblock {\em Concentration inequalities: A nonasymptotic theory of
  independence}.
\newblock Oxford university press.

\bibitem[Bubeck and Cesa-Bianchi, 2012]{bubeck2012regret}
Bubeck, S. and Cesa-Bianchi, N. (2012).
\newblock Regret analysis of stochastic and nonstochastic multi-armed bandit
  problems.
\newblock {\em Foundations and Trends in Machine Learning}, 5(1):1--122.

\bibitem[Bubeck et~al., 2017]{bubeck2017kernel}
Bubeck, S., Lee, Y.~T., and Eldan, R. (2017).
\newblock Kernel-based methods for bandit convex optimization.
\newblock In {\em Proc. 49th Annual ACM SIGACT Symposium on Theory of
  Computing}, pages 72--85.

\bibitem[Duchi et~al., 2015]{duchi2015}
Duchi, J.~C., Jordan, M.~I., Wainwright, M.~J., and Wibisono, A. (2015).
\newblock Optimal rates for zero-order convex optimization: The power of two
  function evaluations.
\newblock {\em IEEE Transactions on Information Theory}, 61(5):2788--2806.

\bibitem[Evans and Gariepy, 2018]{evans2018measure}
Evans, L.~C. and Gariepy, R.~F. (2018).
\newblock {\em Measure theory and fine properties of functions}.
\newblock Routledge.

\bibitem[Flaxman et~al., 2005]{flaxman2004}
Flaxman, A.~D., Kalai, A.~T., and McMahan, H.~B. (2005).
\newblock Online convex optimization in the bandit setting: gradient descent
  without a gradient.
\newblock In {\em Proc. 16th Annual ACM-SIAM Symposium on Discrete algorithms
  (SODA)}, pages 385–--394.

\bibitem[Gasnikov et~al., 2017]{Gasnikov2017}
Gasnikov, A., Krymova, E., Lagunovskaya, A., Usmanova, I., and Fedorenko, F.
  (2017).
\newblock Stochastic online optimization. {S}ingle-point and multi-point
  non-linear multi-armed bandits. {C}onvex and strongly-convex case.
\newblock {\em Automation and Remote Control}, 78(2):224--234.

\bibitem[Gasnikov et~al., 2016]{gasnikov_gradient-free_2016}
Gasnikov, A., Lagunovskaya, A., Usmanova, I., and Fedorenko, F. (2016).
\newblock Gradient-free proximal methods with inexact oracle for convex
  stochastic nonsmooth optimization problems on the simplex.
\newblock {\em Automation and Remote Control}, 77(11):2018--2034.

\bibitem[Ghadimi and Lan, 2013]{Ghadimi2013}
Ghadimi, S. and Lan, G. (2013).
\newblock Stochastic first- and zeroth-order methods for nonconvex stochastic
  programming.
\newblock {\em SIAM Journal on Optimization}, 23(4):2341–2368.

\bibitem[Hu et~al., 2016]{liang2014zeroth}
Hu, X., Prashanth, L.~A., Gy{\"o}rgy, A., and Szepesv{\'a}ri, C. (2016).
\newblock Escaping the local minima via simulated annealing: Optimization of
  approximately convex functions.
\newblock In {\em Proc. 10th International Conference on Artificial
  Intelligence and Statistics}, pages 819--828.

\bibitem[Lattimore and Gyorgy, 2021]{lattimore21a}
Lattimore, T. and Gyorgy, A. (2021).
\newblock Improved regret for zeroth-order stochastic convex bandits.
\newblock In {\em Proceedings of Thirty Fourth Conference on Learning Theory},
  volume 134, pages 2938--2964. PMLR.

\bibitem[Nemirovsky and Yudin, 1983]{NY1983}
Nemirovsky, A.~S. and Yudin, D.~B. (1983).
\newblock {\em Problem Complexity and Method Efficiency in Optimization.}
\newblock Wiley \& Sons.

\bibitem[Nesterov, 2011]{Nesterov2011}
Nesterov, Y. (2011).
\newblock Random gradient-free minimization of convex functions.
\newblock Technical Report 2011001, Center for Operations Research and
  Econometrics (CORE), Catholic University of Louvain.

\bibitem[Nesterov and Spokoiny, 2017]{NS17}
Nesterov, Y. and Spokoiny, V. (2017).
\newblock Random gradient-free minimization of convex functions.
\newblock {\em Found. Comput. Math.}, 17:527--–566.

\bibitem[Novitskii and Gasnikov, 2021]{Gasnikov}
Novitskii, V. and Gasnikov, A. (2021).
\newblock Improved exploiting higher order smoothness in derivative-free
  optimization and continuous bandit.
\newblock {\em arXiv preprint arXiv:2101.03821}.

\bibitem[Orabona, 2019]{Orabona2019}
Orabona, F. (2019).
\newblock A modern introduction to online learning.
\newblock {\em ArXiv}, abs/1912.13213.

\bibitem[Orabona and Pál, 2016]{orabonapal16}
Orabona, F. and Pál, D. (2016).
\newblock Scale-free online learning.
\newblock {\em Theoretical Computer Science}, 716.

\bibitem[Polyak and Tsybakov, 1990]{PT90}
Polyak, T.~B. and Tsybakov, A.~B. (1990).
\newblock Optimal order of accuracy of search algorithms in stochastic
  optimization.
\newblock {\em Problems of Information Transmission}, 26(2):45--53.

\bibitem[Saha and Tewari, 2011]{saha2011}
Saha, A. and Tewari, A. (2011).
\newblock Improved regret guarantees for online smooth convex optimization with
  bandit feedback.
\newblock In {\em Proc. 14th International Conference on Artificial
  Intelligence and Statistics}, pages 636--642.

\bibitem[Schechtman and Zinn, 1990]{Schechtman_Zinn90}
Schechtman, G. and Zinn, J. (1990).
\newblock On the volume of the intersection of two {$L^n_p$} balls.
\newblock {\em Proc. Amer. Math. Soc.}, 110(1):217--224.

\bibitem[Shalev-Shwartz, 2012]{Shalev-Shwartz}
Shalev-Shwartz, S. (2012).
\newblock {\em Online learning and online convex optimization}, volume~4.
\newblock Now Publishers, Inc.

\bibitem[Shamir, 2013]{Shamir13}
Shamir, O. (2013).
\newblock On the complexity of bandit and derivative-free stochastic convex
  optimization.
\newblock In {\em Proc. 30th Annual Conference on Learning Theory}, pages
  1--22.

\bibitem[Shamir, 2017]{Shamir17}
Shamir, O. (2017).
\newblock An optimal algorithm for bandit and zero-order convex optimization
  with two-point feedback.
\newblock {\em Journal of Machine Learning Research}, 18(1):1703--1713.

\bibitem[Zorich, 2016]{zorich2016}
Zorich, V.~A. (2016).
\newblock {\em Mathematical analysis II}.
\newblock Springer.

\end{thebibliography}
\bibliographystyle{apalike}

\appendix

\newpage
\appendix
\begin{center}
    {\Large \bf Supplementary Material}
\end{center}
\vspace{.3truecm}

This supplementary material contains the proofs and results omitted from the main body. In Appendix~\ref{app:IPP} we recall the appropriate version of the Stokes' theorem and discuss its applicability for Lipschitz functions on $\ball_1^d$. In Appendix~\ref{app:variance} we provide the proof of Lemma~\ref{lem:poincare_like_2}. Finally, in Appendix~\ref{app:upper_bound} we provide the proofs of Theorems~\ref{no-noise-reg},~\ref{noisy-reg},~\ref{no-noise-ad},~\ref{noisy-reg-ad}.

\paragraph{Additional notation}

For two functions $g, \eta : \bbR^d \to \bbR$, we denote by $\eta \star g$ their convolution defined point-wise for $\bx \in \bbR^d$ as
\begin{align*}
    \big(\eta \star g\big) (\bx) = \int_{\bbR^d} \eta(\bx - \bx')g(\bx') \d \bx'\enspace.
\end{align*}
The standard mollifier $\eta_{\epsilon}: \bbR^d \to \bbR$ is defined as $\eta_{\epsilon}(\bx) = \epsilon^{-d}\eta_1(\bx / \epsilon)$ for $\epsilon > 0$ and $\bx \in \bbR$, where $\eta_1 : \bbR^d \to \bbR$ is defined as
\begin{align*}
    \eta_1(\bx) = 
    \begin{cases}
    C\exp\parent{\frac{1}{\|\bx\|_2^2-1}} &\text{if } \|\bx\|_2 \leq 1\\
    0 &\text{otherwise}
    \end{cases}\enspace,
\end{align*}
with $C$ chosen so that $\int_{\bbR^d} \eta_1(\bx) \d\bx = 1$.

\section{Integration by parts}
\label{app:IPP}
We first recall the following result that can be found in~\cite[Section 13.3.5, Exercise 14a]{zorich2016}. 

\begin{theorem}[Integration by parts in a multiple integral]
\label{thm:ipp}
    Let $D$ be an open connected subset of $\bbR^d$ with a piecewise smooth boundary $\partial D$ oriented by the outward unit normal $\bn = (n_1, \ldots, n_d)^\top$. Let $g$ be a continuously differentiable  function in $D \cup \partial D$. Then
    \begin{align*}
        \int_{D} \nabla g(\bu) \d \bu = \int_{\partial D} g(\bzeta) \bn(\bzeta) \d S(\bzeta)\enspace.
    \end{align*}
\end{theorem}
\begin{remark}
We refer to~\cite[Section 12.3.2, Definitions 4 and 5]{zorich2016} for the definition of piecewise smooth surfaces and their orientations respectively.
\end{remark}

The idea of using the instance of Theorem \ref{thm:ipp} (also called Stokes' theorem) with $D=\ball^d_2$ to obtain $\ell_2$-randomized estimators of the gradient belongs to~\citet{NY1983}. It was further used in several papers  ~\citep{flaxman2004,BP2016,Shalev-Shwartz,Shamir17} to mention just a few. Those papers were referring to~\cite{NY1983} but~\cite{NY1983} did not provide an exact statement of the result (nor a reference) and only tossed the idea in a discussion. However, the classical analysis formulation
as presented in Theorem \ref{thm:ipp} does not apply to Lipschitz continuous functions that were considered in~\citep{flaxman2004,BP2016,Shalev-Shwartz,Shamir17}. We are not aware of whether its extension to Lipschitz continuous functions, though rather standard, is proved in the literature.

In this paper, we apply Theorem \ref{thm:ipp} with the $\ell_1$-ball $D= \ball^d_1$. Our aim in this section is to provide a variant of Theorem \ref{thm:ipp} applicable to a Lipschitz continuous function $g : \bbR^d \to \bbR$, which is not necessarily continuously differentiable on $D \cup \partial D = \ball^d_1 \cup \sphere^d_1$. To this end, we will go through the argument of approximating $g$ by $C^{\infty}(\Omega)$ functions, where $\Omega \subset \bbR^d$ is an open bounded connected subset of $\bbR^d$ such that $D \cup \partial D \subset \Omega$.
Let $g_n = \eta_{{1}/n} \star g$, where $\eta_{1/n}$ is the standard mollifier. {Let $g : \bbR^d \to \bbR$ be a function satisfying the Lipschitz condition w.r.t. the $\ell_1$-norm: $|g(\bu)-g(\bu')|\le L\|\bu-\bu'\|_1$.}
Since $g$ is continuous in $\Omega$ and, by construction $D \cup \partial D \subset \Omega$, then using basic properties of mollification~\cite[see e.g.,][Theorem 4.1 (ii)]{evans2018measure} we have
\begin{align*}
    g_n \longrightarrow g
\end{align*}
uniformly on $D \cup \partial D$ (in particular, uniformly on $\partial D$).
Furthermore, let $\nabla g$ be the gradient of $g$, which by Rademacher theorem {\cite[see e.g.,][Theorem 3.2]{evans2018measure}} is well defined almost everywhere w.r.t. the Lebesgue measure and
\begin{align*}
    \|\nabla g(\bu)\|_{\infty} \leq L \qquad \text{a.e.}
\end{align*}
It follows that  $\frac{\partial g}{\partial u_j}$ is absolutely integrable on $
\Omega$ for any $j \in [d]$.
Furthermore, since
\begin{align*}
    \frac{\partial g_{n}}{\partial u_j} = \eta_{1/n} \star \left(\frac{\partial g}{\partial u_j}\right)\enspace,
\end{align*}
we can apply~\cite[Theorem 4.1 (iii)]{evans2018measure} that yields
\begin{align*}
    \int_{D}\|\nabla g_{n}(\bu) - \nabla g(\bu)\|_2 \d \bu \longrightarrow 0\enspace.
\end{align*}
{Combining the above remarks we obtain that the result of Theorem \ref{thm:ipp} is valid for functions $g$ that are Lipschitz continuous w.r.t. the $\ell_1$-norm. Thus, it is also valid when the Lipschitz condition is imposed w.r.t. any  $\ell_q$-norm with $q\in[1,\infty]$. Specifying this conclusion for the particular case $D=B_1^d$, we obtain the following theorem.}

\begin{theorem}
\label{thm:ipp2}
Let the function $g: \bbR^d \to \bbR$ {be Lipschitz continuous w.r.t. the $\ell_q$-norm with $q\in[1,\infty]$.} Then
    \begin{align*}
        \int_{\ball^d_1} \nabla g(\bu) \d \bu = \frac{1}{\sqrt{d}}\int_{\sphere^d_1} g(\bzeta) \sign(\bzeta) \d S(\bzeta)\enspace,
    \end{align*}
    where $\nabla g(\cdot)$ is defined up to a set of zero Lebesgue measure by the Rademacher theorem.
\end{theorem}
\section{Proof of Lemma~\ref{lem:poincare_like_2}}
\label{app:variance}

To prove Lemma~\ref{lem:poincare_like_2}, we first recall the weighted Poincaré inequality for the univariate exponential measure (mean $0$ and scale parameter $1$ Laplace distribution).

\begin{lemma}[Lemma 2.1 in~\cite{bobkov1997poincare}]
\label{lem:talagrand}
Let $W$ be mean $0$ and scale parameter $1$ Laplace random variable.
Let $g: \bbR \to \bbR$ be continuous almost everywhere differentiable function such that
\begin{align*}
    \Exp[|g(W)|] < \infty \quad\text{and}\quad \Exp[|g'(W)|] < \infty \quad\text{and}\quad \lim_{|w| \rightarrow \infty} g(w)\exp(-|w|) = 0\enspace,
\end{align*}
then,
\begin{align*}
    \Exp[(g(W) - \Exp[g(W)])^2] \leq 4 \Exp [(g'(W))^2].
\end{align*}
\end{lemma}

We are now in a position to prove Lemma~\ref{lem:poincare_like_2}. The proof is inspired by~\cite[Lemma 2]{barthe2009remarks}.
\begin{proof}[Proof of Lemma~\ref{lem:poincare_like_2}]
Throughout the proof, we assume without loss of generality that $\Exp[G(\bzeta)] = 0$. Indeed, if it is not the case, we use the result for the centered function $\tilde G(\bzeta)= G(\bzeta)-\Exp[G(\bzeta)]$, which has the same gradient. 

First, consider the case of continuously differentiable $G$. 
Let $\bW = \left(W_1, \dots, W_d\right)$ be a vector of i.i.d. mean $0$ and scale parameter $1$ Laplace random variables and define $\bT(\bw) = \bw / \norm{\bw}_1$. Introduce the notation
\begin{align*}
    F(\bw) \eqdef \|\bw\|_1^{\sfrac 1 2}G(\bT(\bw))\enspace.
\end{align*}
Lemma~1 in~\cite{Schechtman_Zinn90} asserts that, for $\bzeta$ uniformly distributed on $\partial B_1^d$,
\begin{align}
    \label{eq:SZ90}
    \bT(\bW) \stackrel{d}{=} \bzeta\quad\text{and}\quad \bT(\bW)  
    \text{ is independent of } \|\bW\|_1 \enspace.
\end{align}
In particular,
\begin{align*}
    \Var(F(\bW)) = d\Var(G(\bzeta))\enspace.
\end{align*}
Using the Efron-Stein inequality~\cite[see e.g.,][Theorem~3.1]{boucheron2013concentration} we obtain
\begin{align*}
    \Var(F(\bW)) \leq \sum_{i = 1}^d\Exp\left[ \Var_{i}(F) \right]\enspace,
\end{align*}
where 
\begin{align*}
    \Var_{i}(F) = \Exp\bigg[ \big(F(\bW) - \Exp[F(\bW) \mid \bW^{-i}] \big)^2 \mid \bW^{-i} \bigg]
\end{align*}
with $\bW^{-i} \eqdef (W_1, \ldots, W_{i - 1}, W_{i+1}, \ldots, W_{d})$. 
Note that on the event $\{\bW^{-i} \neq \bzero\}$ (whose complement has zero measure), the function
\begin{align*}
    w \mapsto F(W_1, \ldots, W_{i - 1}, w, W_{i + 1}, \ldots, W_d)\enspace,
\end{align*}
satisfies the assumptions of Lemma~\ref{lem:talagrand}. Thus, 
\begin{align}
\label{eq:poincare_new_10}
   d\Var(G(\bzeta))= \Var(F(\bW)) \leq 4\sum_{j = 1}^d\Exp\left[ \left(\frac{\partial F}{\partial {w_j}}(\bW)\right)^2 \right] = 4\Exp\|\nabla F(\bW)\|_2^2\enspace.
\end{align}
In order to compute $\nabla F(\bW)$, we observe that for every $i \neq j \in [d]$ we have for all $\bw \neq \bzero$ such that $w_i, w_j \neq 0$
\begin{align*}
    \frac{\partial T_i}{\partial w_j}(\bw) = -\frac{w_i\sign(w_j)}{\|\bw\|_1^2} \quad\text{and}\quad \frac{\partial T_i}{\partial w_i}(\bw) = \frac{1}{\|\bw\|_1} - \frac{w_i\sign(w_i)}{\|\bw\|_1^2}\enspace.
\end{align*}
Thus, the Jacobi matrix of $\bT(\bw)$ has the form
\begin{align*}
    \jac_{\bT}(\bw)
    =
    \frac{\bfI}{\|\bw\|_1} - \frac{{ \bw (\sign(\bw))^\top}}{\|\bw\|_1^2}
    =
    \frac{1}{\|\bw\|_1}\parent{\bfI - \bT(\bw)\big(\sign(\bw)\big)^{\top}}\enspace.
\end{align*}
It follows that almost surely
\begin{align*}
    \nabla F(\bW) = \frac{1}{2\|\bW\|_1^{\sfrac 1 2}}G(\bT(\bW))\sign(\bW) + \frac{1}{\|\bW\|_1^{\sfrac 1 2}}\bigg(\bfI - \bT(\bW)\big(\sign(\bW)\big)^{\top}\bigg)\nabla G(\bT(\bW))\enspace\enspace.
\end{align*}
Observe that since $\scalar{\sign(\bW)}{\bT(\bW)} = 1$ almost surely, we have
\begin{align*}
    \big(\sign(\bW)\big)^{\top}\bigg(\bfI - \bT(\bW)\big(\sign(\bW)\big)^{\top}\bigg)\nabla G(\bT(\bW)) =  0 \quad\text{almost surely}\enspace.
\end{align*}
The above two equations imply that almost surely
\begin{align*}
    4\|\nabla F(\bW)\|_2^2
    &=
    \frac{d}{\|\bW\|_1}G^2(\bT(\bW)) + \frac{4}{\|\bW\|_1}\norm{\bigg(\bfI - \bT(\bW)\big(\sign(\bW)\big)^{\top}\bigg)\nabla G(\bT(\bW))}_2^2\\
    &\leq
    \frac{d}{\|\bW\|_1}G^2(\bT(\bW)) + \frac{4}{\|\bW\|_1}\norm{\nabla G(\bT(\bW))}_2^2(1 + \sqrt{d}\|\bT(\bW)\|_2)^2\enspace,
\end{align*}
where we used the fact that the operator norm of $\bfI - \ba \bb^\top$ is not greater than $1 + \|\ba\|_2 \|\bb\|_2$.
Combining the above bound with~\eqref{eq:poincare_new_10}, and using the facts that $\Exp[\|\bW\|_1^{-1}] = \frac{1}{d-1}$,
 $\Exp[G(\bT(\bW))]=\Exp[G(\bzeta)]=0$ 
and the independence of $\|\bW\|_1$ and $\bT(\bW)$ (cf. \eqref{eq:SZ90}) yields
\begin{align*}
    d\parent{1 - \frac{1}{d-1}}\Var(G(\bzeta)) \leq \frac{4}{d-1}\Exp\left[\|\nabla G(\bT(\bW))\|_2^2{(1 + \sqrt{d}\|\bT(\bW)\|_2)^2}\right]\enspace.
\end{align*}
Rearranging, we deduce the first claim of the lemma since $\bT(\bW) \stackrel{d}{=} \bzeta$.

To prove the second statement of the lemma regarding Lipschitz functions, it is sufficient to apply the first one to $G_{n}$---the sequence of smoothed versions of $G$ such that $G_n \in C^{\infty}(\bbR)$ and
\begin{align*}
    G_{n} \longrightarrow G\enspace,
\end{align*}
uniformly on every compact subset, and $\sup_{n \geq 1}\|\nabla G_n(\bx)\|_2 \leq L$ for almost all $\bx \in \bbR^d$. A sequence $G_n$ satisfying these properties can be constructed by standard mollification due to the fact that $G$ is Lipschitz continuous~\cite[see e.g.,][Theorem 4.2]{evans2018measure}. Finally, to obtain the value $\Exp\|\bT(\bW)\|_2^2 = \Exp\|\bzeta\|_2^2$ we use Lemma~\ref{lem:2ndmoment_l1} below.
\end{proof}

\begin{lemma}
\label{lem:2ndmoment_l1}
Let $\bzeta$ be distributed uniformly on $\sphere^d_1$. Then, $\Exp \norm{\bzeta}^2_2 = \frac{2}{d + 1}$.
\end{lemma}
\begin{proof}
We use the same tools as in the proof of Lemma \ref{eq:moments_uniform_l1}. Let $\bW = (W_1, \ldots, W_d)$ be a vector of \iid random variables following the Laplace distribution with mean $0$ and scale parameter $1$. By \eqref{eq:SZ90} we have that  $\bzeta\stackrel{d}{=} \frac{\bW}{\norm{\bW}_1}$ and $\bzeta$ is independent of $\norm{\bW}_1$. 
Therefore, 
\begin{align}
    \label{eq:2ndmoment_l1_1}
    \Exp\|\bzeta\|_2^2
    = \frac{\Exp\norm{\bW}^2_2}{\Exp\norm{\bW}_1^2}\enspace.
\end{align}
Here,
\begin{align}
\label{eq:2ndmoment_l1_2}
    \Exp\norm{\bW}^2_2 =  \sum_{j = 1}^d \Exp[W_j^2] = d \Exp[W_1^2] = 2d\enspace.
\end{align}
Furthermore, $\norm{\bW}_1$ follows the Erlang distribution with parameters $(d, 1)$, which implies
\begin{align}
\label{eq:2ndmoment_l1_3}
    \Exp\norm{\bW}_1^2 = \frac{1}{\Gamma(d)}\int_{0}^{\infty}x^{d + 1}\exp(-x) \d x =  \frac{\Gamma(d + 2)}{\Gamma(d)}\enspace.
\end{align}
The lemma follows by combining~\eqref{eq:2ndmoment_l1_1} -- \eqref{eq:2ndmoment_l1_3}.
\end{proof}


\section{Upper bounds}\label{app:upper_bound}

The proofs of Theorems~\ref{no-noise-reg},~\ref{noisy-reg},~\ref{no-noise-ad},~\ref{noisy-reg-ad} resemble each other. They only differ in the ways of handling the variance terms depending on $\|\bg_t\|_{p^*}^2$ and in the choice of parameters. For this reason, we suggest the interested reader to follow the proofs in a linear manner starting from the next paragraph. 

\paragraph{Common part of the proofs of Theorems~\ref{no-noise-reg},~\ref{noisy-reg}.}
We start with the part of the proofs that is common for Theorems~\ref{no-noise-reg},~\ref{noisy-reg}. Fix some $\bx \in \com$. 
Due to Assumption~\ref{main-ass}, we can use Lemma~\ref{unb-est}, which implies
\begin{align*}
    \Exp\left[\sum_{t=1}^{T}\scalar{\Exp\left[\bg_t \mid \bx_t\right]}{\bx_t - \bx}\right] 
    = \Exp\left[\sum_{t=1}^{T}\scalar{\nabla {\sf f}_{t, h}(\bx_t)}{\bx_t - \bx}\right]
    \geq \Exp\left[\sum_{t=1}^{T}\big({\sf f}_{t, h}(\bx_t) - {\sf f}_{t, h}(\bx)\big)\right],
\end{align*}
where ${\sf f}_{t, h}(\bx)=\Exp [f_t(\bx +h\bU)]$ with $\bU$ uniformly distributed on $B_1^d$.
Furthermore,
by the approximation property derived in Lemma~\ref{unb-est} and the standard bound on the cumulative regret of dual averaging algorithm~\cite[see e.g.,][Corollary 7.9.]{Orabona2019} we deduce that

\begin{equation}
    \label{eq:meta_main}
\begin{aligned}
    \Exp\left[\sum_{t=1}^{T}\big(f_{t}(\bx_t) - f_{t}(\bx)\big)\right] &\leq \Exp\left[\sum_{t=1}^{T}\langle\Exp\left[\bg_t|\bx_t\right],\bx_t - \bx\rangle\right] + L\mathrm{b}_{q}(d)\sum_{t=1}^{T}h_t
    \\&\leq \frac{R^2}{\eta} + \frac{\eta}{2}\sum_{t=1}^{T}\Exp\norm{\bg_t}_{p^*}^2 + L\mathrm{b}_{q}(d)\sum_{t=1}^{T}h_t\enspace,
\end{aligned}
\end{equation}
where in the last inequality we used the identity $\eta_1 = \ldots = \eta_T = \eta$. The results of Theorems~\ref{no-noise-reg},~\ref{noisy-reg} follow from the bound~\eqref{eq:meta_main} as detailed below.

\begin{proof}[Proof of Theorem~\ref{no-noise-reg}]
Here $h_1 = \ldots = h_T = h$, and we work under Assumption~\ref{ass:no_noise}.
In this case, bounding $\Exp\|\bg_t\|_{p^*}$ in ~\eqref{eq:meta_main} via Lemma~\ref{var} yields
\begin{align*}
    \Exp\left[\sum_{t=1}^{T}\big(f_{t}(\bx_t) - f_{t}(\bx)\big)\right] &\leq {\frac{R^2}{\eta} + 6(1+\sqrt{2})^2L^2\cdot\eta T d^{1+\frac{2}{q \wedge 2}-\frac{2}{p}}}+ LhT\mathrm{b}_{q}(d)\enspace.
\end{align*}
Minimizing the the right hand side of the above inequality over $\eta > 0$ and substituting $\eta = \frac{R}{L\left(\sqrt{6}+\sqrt{12}\right)}\sqrt{\frac{d^{-1-\frac{2}{{q \wedge 2}}+\frac{2}{p}}}{T}}$ we deduce that
\begin{align*}
    \Exp\left[\sum_{t=1}^{T}\big(f_{t}(\bx_t) - f_{t}(\bx)\big)\right] &\leq  2\left(\sqrt{6} + \sqrt{12}\right)RL
    d^{\frac{1}{2}+\frac{1}{{q \wedge 2}}-\frac{1}{p}}\sqrt{T}+LhT\mathrm{b}_{q}(d)\enspace.
\end{align*}
Taking $h \leq \frac{7R}{100\mathrm{b}_{q}(d)\sqrt{T}}d^{\frac{1}{2}+\frac{1}{q \wedge 2}-\frac{1}{p}}$ makes negligible the second summand in the above bound. This concludes the proof.
\end{proof}

\begin{proof}[Proof of Theorem~\ref{noisy-reg}]
Here again $h_1 = \ldots = h_T = h$, but we work under Assumption~\ref{ass:noise}.
Then, bounding $\Exp\|\bg_t\|_{p^*}$ in~\eqref{eq:meta_main} via Lemma~\ref{var} yields
\begin{align*}
    \Exp\left[\sum_{t=1}^{T}\big(f_{t}(\bx_t) - f_{t}(\bx)\big)\right] \leq {\frac{R^2}{\eta} + 
    \eta T\left(\frac{d^{4-\frac{2}{p}}\sigma^2}{h^2}+6\left(1+\sqrt{2}\right)^2L^2d^{1+\frac{2}{{q \wedge 2}}-\frac{2}{p}}\right)} + LhT\mathrm{b}_{q}(d)\enspace.
\end{align*}
Minimizing the right hand side of the above inequality over $\eta > 0$ and substituting the optimal value
\[ \eta = \frac{R}{\sqrt{T}}\left(\frac{d^{4-\frac{2}{p}}\sigma^2}{2h^2}+6\left(1+\sqrt{2}\right)^2L^2d^{1+\frac{2}{{q \wedge 2}}-\frac{2}{p}}\right)^{-\frac{1}{2}}\enspace,\] results in the following upper bound on the regret
\begin{equation*}
\begin{aligned}
    \Exp\left[\sum_{t=1}^{T}\big(f_{t}(\bx_t) - f_{t}(\bx)\big)\right]
    &\leq
    2R\sqrt{T}\left(\frac{d^{4-\frac{2}{p}}\sigma^2}{2h^2}+6\left(1+\sqrt{2}\right)^2L^2d^{1+\frac{2}{{q \wedge 2}}-\frac{2}{p}}\right)^{\frac{1}{2}} + LhT\mathrm{b}_{q}(d)
    \\
    &\leq
    2\left(\sqrt{6}+\sqrt{12}\right)RL\sqrt{Td^{1+\frac{2}{{q \wedge 2}}-\frac{2}{p}}}+{\sqrt{2}R\sqrt{T}\frac{d^{2-\frac{1}{p}}\sigma}{h} + LhT\mathrm{b}_{q}(d)}\enspace,
\end{aligned}
\end{equation*}
where for the last inequality we used the fact that $\sqrt{a + b} \leq \sqrt{a} + \sqrt{b}$ for $a, b \geq 0$.
Minimizing over $h > 0$ the last expression and substituting the optimal value $h = \left(\frac{\sqrt{2}R\sigma}{L\mathrm{b}_{q}(d)}\right)^{\frac{1}{2}}T^{-\frac{1}{4}}d^{1-\frac{1}{2p}}$ we get
\begin{equation*}
\begin{aligned}
    \Exp\left[\sum_{t=1}^{T}\big(f_{t}(\bx_t) - f_{t}(\bx)\big)\right] \leq 11.9 RL\sqrt{Td^{1+\frac{2}{{q \wedge 2}}-\frac{2}{p}}}+2.4\sqrt{RL\sigma}T^{\frac{3}{4}} \sqrt{\mathrm{b}_{q}(d)}d^{\frac{1}{2}-\frac{1}{2p}}.\qedhere
\end{aligned}
\end{equation*}
\end{proof}
\paragraph{Common part of the proofs of Theorems~\ref{no-noise-ad},~\ref{noisy-reg-ad}.} Here, we state the common parts of the proofs for Theorems~\ref{no-noise-ad},~\ref{noisy-reg-ad}. Similar to the first inequality in \eqref{eq:meta_main}, we have
\begin{align*}
    \Exp\left[\sum_{t=1}^{T}\big(f_{t}(\bx_t) - f_{t}(\bx)\big)\right] &\leq \Exp\left[\sum_{t=1}^{T}\langle\bg_t,\bx_t - \bx\rangle\right] + L\mathrm{b}_{q}(d)\sum_{t=1}^{T}h_t\enspace.
\end{align*}
Note that without loss of generality, we can assume that $\sum_{k=1}^{t}\norm{\bg_k}_{p^*}^{2} \neq 0$, for all $t\geq1$. This is a consequence of the fact that if $\sum_{k=1}^{t}\norm{\bg_k}_{p^*}^{2} = 0$, then the first term on the r.h.s. of the above inequality will be zero up to round $t$. Thus, we can erase these iterates from the cumulative regret, only paying the bias term for those rounds. In what follows we essentially use~\cite[Corollary 1]{orabonapal16}, which we re-derive for the sake of clarity. Assume that $\eta_t = \tfrac{\lambda}{\sqrt{\sum_{k=1}^{t-1}\norm{\bg_k}^{2}_{p^*}}}$ for $t\in\{2,\dots,T\}$ and $\lambda >0$. Then, applying \cite[Theorem 1]{orabonapal16} we deduce that
\begin{align*}
    \Exp\left[\sum_{t=1}^{T}\big(f_{t}(\bx_t) - f_{t}(\bx)\big)\right] &\leq \left(\frac{R^2}{\lambda} + 2.75\cdot\lambda\right)\Exp\left[\sqrt{\sum_{t=1}^{T}\norm{\bg_t}^2_{p^*}}\right] \\&\phantom{\leq}+ 3.5D\cdot \Exp[\max_{t\in[T]}\norm{\bg_t}_{p^*}]+ L\mathrm{b}_{q}(d)\sum_{t=1}^{T}h_t\enspace,
\end{align*}
where we introduced $D = \sup_{\bu,\bw\in\com}\norm{\bu-\bw}_{p}$. By \cite[Proposition 1]{orabonapal16}, we have $D\leq \sqrt{8}R$. Moreover, by Jensen's inequality, using the rough bound $\Exp[\max_{t\in[T]}\norm{\bg_t}_{p^*}]\leq \sqrt{\sum_{t=1}^{T}\Exp\left[\norm{\bg_t}^2_{p^*}\right]}$, and substituting $\lambda = \frac{R}{\sqrt{2.75}}$, we deduce that 
\begin{align}\label{eq:meta_ad}
    \Exp\left[\sum_{t=1}^{T}\big(f_{t}(\bx_t) - f_{t}(\bx)\big)\right] &\leq \left(2\sqrt{2.75} + 3.5\sqrt{8}\right)R\sqrt{\sum_{t=1}^{T}\Exp\left[\norm{\bg_t}^{2}_{p^*}\right]} +L\mathrm{b}_{q}(d)\sum_{t=1}^{T}h_t\enspace.
\end{align}
Proofs of Theorems~\ref{no-noise-ad},~\ref{noisy-reg-ad} provided below follow from the above inequality by properly selecting $h_t > 0$.
\begin{proof}[Proof of Theorem~\ref{no-noise-ad}]
The bound of Lemma~\ref{var} under Assumption~\ref{ass:no_noise} applied to \eqref{eq:meta_ad} yields
\begin{align*}
    \Exp\left[\sum_{t=1}^{T}\big(f_{t}(\bx_t) - f_{t}(\bx)\big)\right] &\leq 2\left(2\sqrt{2.75} + 3.5\sqrt{8}\right)\left(\sqrt{3}+\sqrt{6}\right)RL\sqrt{T d^{1+\frac{2}{q\wedge 2}-\frac{2}{p}}}  + L\mathrm{b}_{q}(d)\sum_{t=1}^{T}h_t
    \\&\leq 110.53\cdot RL\sqrt{T d^{1+\frac{2}{q\wedge 2}-\frac{2}{p}}}  + L\mathrm{b}_{q}(d)\sum_{t=1}^{T}h_t\enspace.
\end{align*}
Taking $h_t \leq \frac{7R}{200\mathrm{b}_{q}(d)\sqrt{t}}d^{\frac{1}{2}+\frac{1}{q \wedge 2}-\frac{1}{p}}$ makes negligible the last summand in the above bound. This concludes the proof.
\end{proof}

\begin{proof}[Proof of Theorem~\ref{noisy-reg-ad}]
Using~\eqref{eq:meta_ad}, the bound of Lemma~\ref{var} under Assumption~\ref{ass:noise} and the fact that $\sqrt{a + b} \leq \sqrt{a} + \sqrt{b}$ for $a, b \geq 0$, we deduce that
\begin{align*}
    \Exp\left[\sum_{t=1}^{T}\big(f_{t}(\bx_t) - f_{t}(\bx)\big)\right]&\leq  \left(2\sqrt{2.75}+3.5\sqrt{8}\right)R\left(\sum_{t=1}^{T}\frac{d^{4-\frac{2}{p}}\sigma^2}{h_t^2} + 12(1+\sqrt{2})^2L^2T \cdot d^{1 + \frac{2}{q \wedge 2} - \frac{2}{p}}\right)^{\frac{1}{2}}\\
    &\phantom{\leq}+ L\mathrm{b}_{q}(d)\sum_{t=1}^{T}h_t
    \\&\leq 110.6\cdot RL\sqrt{T d^{1+\frac{2}{q\wedge 2}-\frac{2}{p}}} + 13.3R\cdot d^{2-\frac{1}{p}}\sigma \left(\sum_{t=1}^{T} \frac{1}{h_t^{2}}\right)^\frac{1}{2}\\
    &\phantom{\leq}+  L\mathrm{b}_{q}(d)\sum_{t=1}^{T}h_t\enspace.
\end{align*}
Since $h_t = \left(6.65\sqrt{6}\cdot\frac{R}{\mathrm{b}_{q}(d)}\right)^{\frac{1}{2}}t^{-\frac{1}{4}}d^{1-\frac{1}{2p}}$ and $\sum_{t=1}^T t^{\frac{1}{2}} \leq \tfrac{2}{3} T^{\frac{3}{2}}$ and $\sum_{t=1}^{T}t^{-\frac{1}{4}}\leq \tfrac{4}{3} T^{\frac{3}{4}}$, we get
\begin{align*}
   \Exp\left[\sum_{t=1}^{T}\big(f_{t}(\bx_t) - f_{t}(\bx)\big)\right] 
   &\leq
   110.6\cdot RL\sqrt{T d^{1+\frac{2}{q\wedge 2}-\frac{2}{p}}} +5.9\cdot\sqrt{R}\left(\sigma+L\right)T^{\frac{3}{4}} \sqrt{\mathrm{b}_{q}(d)}d^{\frac{1}{2}-\frac{1}{2p}}\enspace.\qedhere
\end{align*}
\end{proof}

\section{Definition of $\ell_2$-randomized estimator}
\label{app:l2_description}
In this section we recall the algorithm of~\citet{Shamir17}. Let $\bzeta^{\circ} \in \bbR^d$ be distributed uniformly on $\sphere^d_2$. Instead of the gradient estimator that we introduce in Algorithm~\ref{algo:simplex}, at a each step $t \geq 1$,~\citet{Shamir17} uses
\begin{align*}
    \bg_t^{\circ} \triangleq \frac{d}{2h}(y_t' - y_t'')\bzeta_t^{\circ}\enspace,
\end{align*}
where $y_t' = f_t(\bx_t + h_t\bzeta^{\circ})$, $y_t'' = f_t(\bx_t - h_t \bzeta_t^{\circ})$, and $\bzeta^{\circ}_t$'s are  independent random variables with the same distribution as $\bzeta^{\circ}$. 

\end{document}